\documentclass[12pt,leqno,oneside]{amsproc}
\usepackage[T2A]{fontenc}
\usepackage[utf8]{inputenc}
\usepackage[russian]{babel}
\usepackage{amssymb}
\allowdisplaybreaks

\renewcommand{\subsection}{\refstepcounter{subsection}%
\par\bigskip\noindent\textbf{\upshape\arabic{subsection}. }}
\renewcommand{\subsubsection}{\refstepcounter{subsubsection}%
\par\medskip\noindent\textbf{\upshape\arabic{subsection}.%
\arabic{subsubsection}.\ }}
\renewcommand{\paragraph}{\refstepcounter{paragraph}%
\par\smallskip\noindent\textbf{\upshape\arabic{subsection}.%
\arabic{subsubsection}.\arabic{paragraph}.\ }}

\numberwithin{equation}{subsection}

\newcommand{\Ya}[1]{\text{\upshape Я}_{#1}}

\title{О сравнении интегралов Дарбу и Римана в конструктивном математическом
анализе}
\author{А.~А.~Владимиров}

\begin{document}
\renewcommand{\proofname}{\upshape Д\:о\:к\:а\:з\:а\:т\:е\:л\:ь\:с\:т\:в\:о}
\pagestyle{plain}
\begin{abstract}
В статье устанавливается, что эквивалентные с точки зрения "`классической"'
математики определения интеграла функции одного вещественного аргумента
по Дарбу и по Риману в рамках конструктивного математического анализа
А.~А.~Маркова оказываются существенно различными.
\end{abstract}
\begin{flushleft}
УДК~510.25
\end{flushleft}
\maketitle

\section{Введение}\label{pt:1}
\subsection
В рамках "`классической"' теоретико-множественной математики известно два основных
способа определения одномерного интеграла Римана: в качестве предела интегральных
сумм при исчезновении измельчённости дробления (подход Римана), и в качестве
промежуточного значения множеств интегралов "<элементарно интегрируемых"> функций,
оценивающих рассматриваемую (подход Дарбу). С точки зрения "`классической"'
математики указанные два подхода приводят к равнообъёмным понятиям. Основная цель
настоящей статьи заключается в построении примера всюду на отрезке \([0,1]\)
заданной функции, интегрируемой по Риману, но не интегрируемой по Дарбу, с точки
зрения конструктивного математического анализа.

В дальнейшем мы всегда, не оговаривая этого особо, будем исходить
из конструктивного понимания математических суждений, как оно даётся ступенчатой
семантической системой А.~А.~Маркова \cite{Markov:1974}, \cite{VD:2009}.

\subsection\label{pt:1:1}
Конструктивный вариант определения Римана рассматривался неоднократно
и в настоящее время может считаться стандартным \cite[Гл.~7, \S~1]{Ku:1973}.
Конструктивную интегрируемость по Дарбу мы вводим следующим образом:

\subsubsection\label{darb}
{\itshape Функция \(f:[0,1]\to\mathbb R\) называется интегрируемой по Дарбу, если
для любого вещественного числа \(\varepsilon>0\) осуществимы две полигональные
функции \(f^{\downarrow}_{\varepsilon}\) и \(f^{\uparrow}_{\varepsilon}\),
удовлетворяющие соотношениям
\begin{gather*}
	(\forall x\in [0,1])\qquad f^{\downarrow}_{\varepsilon}(x)
	\leqslant f(x)\leqslant f^{\uparrow}_{\varepsilon}(x),\\
	\int\limits_0^1(f^{\uparrow}_{\varepsilon}-
	f^{\downarrow}_{\varepsilon})<\varepsilon.
\end{gather*}
}

Здесь и далее функцию \(f:[0,1]\to\mathbb R\) мы называем \emph{полигональной}
(ср.~\cite[Гл.~5, \S~1, Определение~7]{Ku:1973}, \cite[\S~1]{ZTs:1962}), если
осуществимы список \(\{q_k\}_{k=0}^m\) рациональных чисел и упорядоченный
по возрастанию список \(\{p_k\}_{k=0}^m\) попарно различных рациональных точек
отрезка \([0,1]\), удовлетворяющие соотношениям \(p_0=0\), \(p_m=1\) и
\[
	(\forall k\in\mathbb N:1\leqslant k\leqslant m)\,(\forall x\in
		[p_{k-1},p_k])\qquad f(x)=\dfrac{q_k\cdot (x-p_{k-1})+
		q_{k-1}\cdot (p_k-x)}{p_k-p_{k-1}}.
\]
Под \emph{полигональным интегралом} такой функции, как обычно \cite[Гл.~8, \S~2,
Определение~6]{Ku:1973}, \cite[\S~1]{ZTs:1962}, понимается величина
\[
	\int\limits_0^1 f\rightleftharpoons
		\sum\limits_{k=1}^m\dfrac{(q_k+q_{k-1})\cdot (p_k-p_{k-1})}{2},
\]
численно совпадающая с интегралом Римана.


\section{Вспомогательные утверждения}\label{pt:2}
\subsection
Очевидным образом имеет место следующий факт:

\subsubsection\label{prop:2:1:2}
{\itshape Пусть функция \(f:[0,1]\to\mathbb R\) интегрируема по Дарбу. Тогда
для любого рационального числа \(\varepsilon>0\) осуществимы рациональное число
\(\delta>0\) и неотрицательная полигональная функция \(g\),
удовлетворяющие соотношениям
\begin{gather*}
	\int\limits_0^1 g<1/2,\\
	(\forall x,y\in [0,1] : \sup(g(x),\,g(y))<1,\,|x-y|<\delta)\qquad
		|f(x)-f(y)|<\varepsilon.
\end{gather*}
}

\subsection\label{pt:2:2}
Следующее утверждение представляет собой удобный для наших дальнейших целей
аналог известной \cite[Гл.~8, \S~1, Теорема~2]{Ku:1973}, \cite[Теорема~2.1]{%
ZTs:1962} теоремы о сингулярных покрытиях. Оно может также быть рассмотрено
в качестве варианта конструктивного опровержения \cite[Замечание~2]{Dem:1973}
теоремы Лебега о мажорируемой сходимости.

\subsubsection\label{prop:2:2:1}
{\itshape Осуществима неубывающая последовательность \(\{h_n\}_{n=0}^{\infty}\)
неотрицательных полигональных функций, удовлетворяющая соотношениям
\begin{gather}\label{eq:2:1:20}
	(\forall x\in [0,1])\qquad \lim\limits_{n\to\infty} h_n(x)=2,\\
	\label{eq:2:1:21}
	(\forall n\in\mathbb N)\qquad \int\limits_0^1 h_n<1/2.
\end{gather}
}

\begin{proof}
Зафиксируем \cite[Гл.~8, \S~1, Теорема~2]{Ku:1973} накрывающую отрезок \([0,1]\)
последовательность \(\{(a_n,b_n)\}_{n=0}^{\infty}\) непустых интервалов
с рациональными концами, удовлетворяющую соотношению
\begin{equation}\label{eq:2:1:1}
	(\forall n\in\mathbb N)\qquad \sum\limits_{k=0}^n (b_k-a_k)<1/6.
\end{equation}
Введём в рассмотрение последовательность \(\{\varphi_n\}_{n=0}^{\infty}\)
полигональных функций вида
\begin{equation}\label{eq:2:1:2}
	(\forall n\in\mathbb N)\,(\forall x\in [0,1])\qquad \varphi_n(x)=
		\inf\left(1,\,\sup\left(0,\,2-\dfrac{|2x-b_n-a_n|}{b_n-a_n}
		\right)\right),
\end{equation}
а также неубывающую последовательность \(\{h_n\}_{n=0}^{\infty}\) неотрицательных
полигональных функций, рекуррентно заданную соотношениями
\begin{align}\label{eq:2:1:3}
	(\forall x\in [0,1])&& h_0(x)&=0,\\ \label{eq:2:1:4}
	(\forall n\in\mathbb N)\,(\forall x\in [0,1])&& h_{n+1}(x)&=
		\sup(h_n(x),\,2\varphi_n(x)).
\end{align}
Тогда при любом выборе индекса \(n\in\mathbb N\) выполняются неравенства
\begin{flalign*}
	&& \int\limits_0^1 h_{n+1}&\leqslant\sum\limits_{k=0}^n
		\int\limits_0^1 2\varphi_k&\text{[\eqref{eq:2:1:3},
		\eqref{eq:2:1:4}]}&\\
	&& &\leqslant\sum\limits_{k=0}^n 3\cdot(b_k-a_k)&
		\text{[\eqref{eq:2:1:2}]}&\\
	&& &<1/2.&\text{[\eqref{eq:2:1:1}]}&
\end{flalign*}
Кроме того, выполняется соотношение
\begin{flalign*}
	&& (\forall n\in\mathbb N)\,(\forall x\in [0,1])\qquad
		h_n(x)&\leqslant 2,&\text{[\eqref{eq:2:1:3}, \eqref{eq:2:1:4},
		\eqref{eq:2:1:2}]}&
\end{flalign*}
а также обусловленное вложением \([0,1]\subseteq\bigcup
\limits_{n=0}^{\infty} (a_n,b_n)\) соотношение
\begin{flalign*}
	&& (\forall x\in [0,1])\,(\exists n\in\mathbb N)\qquad
		h_{n+1}(x)&\geqslant 2.&\text{[\eqref{eq:2:1:4},
		\eqref{eq:2:1:2}]}&
\end{flalign*}
Тем самым, функциональная последовательность \(\{h_n\}_{n=0}^{\infty}\)
удовлетворяет всем предъявленным в формулировке доказываемого утверждения
требованиям.
\end{proof}

\subsection
Произвольно фиксированным полигональной функции \(f\) и вещественному числу
\(\varepsilon>0\) можно сопоставить \cite[\S~11.2.2]{Sh:1962} функцию
\(\omega(f,\varepsilon):[0,1]\to\mathbb R\) вида
\[
	(\forall x\in [0,1])\qquad [\omega(f,\varepsilon)](x)=
		\sup\limits_{t,s\in [x-\varepsilon,x+\varepsilon]\cap [0,1]}
		|f(t)-f(s)|.
\]
Очевидным образом имеют место следующие два факта:

\subsubsection
{\itshape Пусть даны полигональная функция \(f\) и вещественное число
\(\varepsilon>0\). Тогда функция \(\omega(f,\varepsilon):[0,1]\to\mathbb R\)
является равномерно непрерывной.
}

\subsubsection\label{prop:2:2:2}
{\itshape Пусть даны полигональная функция, вещественное число \(\varepsilon>0\),
а также имеющее не превосходящую \(\varepsilon\) измельчённость интегральное
дробление \(\tau\) отрезка \([0,1]\). Тогда отвечающее дроблению \(\tau\) значение
\(I(f,\tau)\) интегральной суммы функции \(f\) удовлетворяет соотношению
\[
	\left|I(f,\tau)-\int\limits_0^1 f\right|
		\leqslant\int\limits_0^1\omega(f,\varepsilon).
\]
}

Также имеют место следующие два факта:

\subsubsection\label{prop:2:2:3}
{\itshape Пусть \(\{f_n\}_{n=0}^{\infty}\) "--- последовательность полигональных
функций, поточечно сходящаяся к некоторой функции \(f:[0,1]\to\mathbb R\). Пусть
также для любого вещественного числа \(\varepsilon>0\) осуществимо вещественное
число \(\delta>0\), удовлетворяющее соотношению
\[
	(\forall n\in\mathbb N)\qquad \int\limits_0^1\omega(f_n,\delta)
		<\varepsilon.
\]
Тогда функция \(f\) интегрируема по Риману.
}

\begin{proof}
Зафиксируем произвольное вещественное число \(\varepsilon>0\), а также вещественное
число \(\delta>0\), удовлетворяющее соотнощению
\begin{equation}\label{eq:2:1:11}
	(\forall n\in\mathbb N)\qquad \int\limits_0^1\omega(f_n,\delta)
		<\varepsilon/3.
\end{equation}
Тогда для любых двух интегральных дроблений \(\tau\) и \(\sigma\) отрезка
\([0,1]\), имеющих не превосходящую \(\delta\) измельчённость, будут выполняться
соотношения
\begin{flalign*}
	&& |I(f,\tau)-I(f,\sigma)|&=\lim\limits_{n\to\infty} |I(f_n,\tau)-
		I(f_n,\sigma)|&\\
	&& &\leqslant 2\varepsilon/3&\text{[\ref{prop:2:2:2},
		\eqref{eq:2:1:11}]}&\\
	&& &<\varepsilon.
\end{flalign*}
Интегрируемость функции \(f\) по Риману вытекает теперь из произвольности выбора
вещественного числа \(\varepsilon>0\) и критерия Коши.
\end{proof}

\subsubsection\label{prop:2:2:4}
{\itshape Пусть \(\alpha>0\), \(\beta>0\) и \(\zeta\in (0,1)\) "--- рациональные
числа, и пусть полигональная функция \(f\) имеет вид
\[
	(\forall x\in [0,1])\qquad f(x)=\alpha\cdot\sup\left(0,\,1-
		\dfrac{|x-\zeta|}{\beta}\right).
\]
Тогда при любом выборе рационального числа \(\varepsilon>0\) выполняется
соотношение
\begin{equation}\label{eq:2:2:1}
	\int\limits_0^1\omega(f,\varepsilon)\leqslant 8\alpha\varepsilon.
\end{equation}
}

\begin{proof}
Заметим, что всегда выполняется одно из неравенств \(\varepsilon<\beta\)
или \(\varepsilon\geqslant\beta\).

В первом случае функция \(f\) является липшицевой с коэффициентом \(\alpha/\beta\),
а потому функция \(\omega(f,\varepsilon)\) мажорируется постоянной
\(2\alpha\varepsilon/\beta\). Кроме того, указанная функция обращается в нуль
вне отрезка \([\zeta-2\beta,\zeta+2\beta]\). Тем самым, неравенство
\eqref{eq:2:2:1} выполняется.

Во втором случае функция \(\omega(f,\varepsilon)\) мажорируется постоянной
\(\alpha\) и обращается в нуль вне отрезка \([\zeta-2\varepsilon,\zeta+
2\varepsilon]\). Тем самым, неравенство \eqref{eq:2:2:1} также выполняется.
\end{proof}


\section{Построение примера}\label{pt:3}
\subsection
На протяжении настоящего параграфа мы будем считать, что натуральные и рациональные
числа представляют собой слова в трёхбуквенном алфавите \(\{{|},{-},{/}\}\)
\cite[\S~1.6]{MN:1996}, а списки рациональных чисел являются
\mbox{\(*\)-сис}\-те\-ма\-ми \cite[\S~24]{MN:1996}. Поскольку любое конструктивное
отображение множества натуральных чисел в множество \mbox{\(*\)-сис}\-тем
рациональных чисел может быть задано посредством нормального алгорифма
в шестибуквенном алфавите \(\{{|},{-},{/},{*},{a},{b}\}\)
\cite[\S~41.7.1]{MN:1996}, то в дальнейшем мы будем ограничиваться рассмотрением
именно таких алгорифмов.

\subsection\label{pt:3:2}
Зафиксируем некоторую нумерацию \(\{\mathfrak A_n\}_{n=0}^{\infty}\) нормальных
алгорифмов указанного в предыдущем пункте вида. Рассмотрим связанное с этой
нумерацией множество индексов \(N\), для которых процесс применения алгорифма
\(\mathfrak A_N\) к слову \(N\) останавливается с результатом вида
\begin{equation}\label{eq:3:1:1}
	{*}\delta{*}p_0{*}q_0{*}\ldots {*}p_m{*}q_m{*},
\end{equation}
где \(\delta\) "--- положительное рациональное число, а списки рациональных чисел
\(\{p_k\}_{k=0}^m\) и \(\{q_k\}_{k=0}^m\) определяют [\ref{pt:1}.\ref{pt:1:1}]
некоторую неотрицательную полигональную функцию \(g\) со свойством
\begin{equation}\label{eq:3:1:2}
	\int\limits_0^1 g<1/2.
\end{equation}
Указанное множество индексов с очевидностью является бесконечным и полуразрешимым,
а потому допускает перечисление без повторений некоторым полным арифметическим
алгорифмом \(\mu:\mathbb N\to\mathbb N\).

Последовательность рациональных чисел \(\delta\) из представления \eqref{eq:3:1:1}
результатов применения алгорифмов \(\mathfrak A_{\mu(n)}\) к словам \(\mu(n)\)
мы на протяжении настоящего параграфа будем обозначать в виде \(\{\delta_n\}_{%
n=0}^{\infty}\). Последовательность соответствующих полигональных функций мы будем
обозначать в виде \(\{g_n\}_{n=0}^{\infty}\). Кроме того, мы будем предполагать
зафиксированной неубывающую последовательность \(\{h_n\}_{n=0}^{\infty}\)
неотрицательных полигональных функций из утверждения \ref{pt:2}.\ref{prop:2:2:1}.

\subsection
Имеют место следующие три факта:

\subsubsection\label{prop:3:3:1}
{\itshape Осуществимы возрастающая последовательность \(\nu:\mathbb N\to
\mathbb N\), а также последовательность \(\{\zeta_n\}_{n=0}^{\infty}\) рациональных
точек интервала \((0,1)\) и последовательность \(\{\beta_n\}_{n=0}^{\infty}\)
положительных рациональных чисел, удовлетворяющие соотношениям
\begin{gather}\label{eq:3:3:1}
	(\forall n\in\mathbb N)\qquad \beta_n<\inf(\delta_n,\,
		\zeta_n,\,1-\zeta_n),\\ \label{eq:3:3:2}
	(\forall n\in\mathbb N)\,(\forall x\in [\zeta_n-\beta_n,\zeta_n+
		\beta_n])\qquad [g_n+h_{\nu(n)}](x)<1<h_{\nu(n+1)}(x).
\end{gather}
}

\begin{proof}
Положим \(\nu(0)\rightleftharpoons 0\). Построение значений \(\zeta_n\),
\(\beta_n\) и \(\nu(n+1)\) на основе известного значения \(\nu(n)\) может теперь
быть произведено следующим образом. В качестве \(\zeta_n\) выберем произвольную
рациональную точку интервала \((0,1)\), удовлетворяющую неравенству \([g_n+
h_{\nu(n)}](\zeta_n)<1\). Осуществимость такой точки гарантирована соотношениями
\begin{flalign*}
	&& \int\limits_0^1 (g_n+h_{\nu(n)})&<1/2+1/2&
		\text{[\ref{pt:3:2}\,\eqref{eq:3:1:2},
		\ref{pt:2}.\ref{pt:2:2}\,\eqref{eq:2:1:21}]}&\\
	&& &=1.
\end{flalign*}
В качестве \(\nu(n+1)\) выберем произвольное натуральное число, удовлетворяющее
неравенству \(h_{\nu(n+1)}(\zeta_n)>1\). Осуществимость такого числа гарантирована
соотношением \ref{pt:2}.\ref{pt:2:2}\,\eqref{eq:2:1:20}. В качестве \(\beta_n\)
теперь остаётся выбрать произвольное положительное рациональное число, чья малость
будет достаточна для выполнения соотношений \eqref{eq:3:3:1} и \eqref{eq:3:3:2}.
\end{proof}

\subsubsection\label{prop:3:3:2}
{\itshape Пусть \(\{f_n\}_{n=0}^{\infty}\) "--- последовательность полигональных
функций вида
\begin{equation}\label{eq:3:3:3}
	(\forall n\in\mathbb N)\,(\forall x\in [0,1])\qquad f_n(x)=
		2^{-\mu(n)}\cdot\sup\left(0,\,1-\dfrac{|x-\zeta_n|}{\beta_n}
		\right),
\end{equation}
где \(\{\zeta_n\}_{n=0}^{\infty}\) и \(\{\beta_n\}_{n=0}^{\infty}\) "---
числовые последовательности из утверждения \ref{prop:3:3:1}. Тогда осуществима
и интегрируема по Риману функция \(f:[0,1]\to\mathbb R\), удовлетворяющая
соотношению
\[
	(\forall x\in [0,1])\qquad f(x)=\sum\limits_{k=0}^{\infty} f_k(x).
\]
}

\begin{proof}
Зафиксируем произвольную точку \(x\in [0,1]\), а также номер \(n\in\mathbb N\),
удовлетворяющий неравенству \(h_{\nu(n)}(x)>1\)
[\ref{pt:2}.\ref{pt:2:2}\,\eqref{eq:2:1:20}]. Тогда выполняется соотношение
\begin{flalign*}
	&& (\forall k\in\mathbb N: k\geqslant n)\qquad f_k(x)&=0,&
		\text{[\eqref{eq:3:3:3}, \eqref{eq:3:3:2}]}&
\end{flalign*}
означающее сходимость ряда
\[
	\sum\limits_{k=0}^{\infty} f_k(x).
\]
Кроме того, для любых вещественного числа \(\varepsilon>0\) и номера
\(n\in\mathbb N\) выполняются оценки
\begin{flalign*}
	&& \int\limits_0^1\omega\left(\sum\limits_{k=0}^n f_k,\,\varepsilon/16
		\right)&\leqslant\sum\limits_{k=0}^n\int\limits_0^1\omega(f_k,\,
		\varepsilon/16)\\
	&& &<\sum\limits_{k=0}^n 2^{-\mu(k)-1}\cdot\varepsilon&
		\text{[\eqref{eq:3:3:3}, \ref{pt:2}.\ref{prop:2:2:4}]}&\\
	&& &<\varepsilon,
\end{flalign*}
означающие интегрируемость функции \(f\) по Риману [\ref{pt:2}.\ref{prop:2:2:3}].
\end{proof}

\subsubsection
{\itshape Функция \(f:[0,1]\to\mathbb R\) из утверждения \ref{prop:3:3:2}
не является интегрируемой по Дарбу.
}

\begin{proof}
В случае интегрируемости рассматриваемой функции по Дарбу должен найтись
[\ref{pt:2}.\ref{prop:2:1:2}] алгорифм, перерабатывающий каждое натуральное число
\(n\) в список вида \ref{pt:3:2}\,\eqref{eq:3:1:1}, отвечающие которому
положительное рациональное число \(\delta>0\) и неотрицательная полигональная
функция \(g\) удовлетворяют соотношению
\[
	(\forall x,y\in [0,1]: \sup(g(x),\,g(y))<1,\,|x-y|<\delta)\qquad
		|f(x)-f(y)|<2^{-n-1}.
\]
При этом, очевидно, найдётся также натуральное число \(m\), для которого значение
\(\mu(m)\) будет являться номером рассматриваемого алгорифма при нумерации
из пункта \ref{pt:3:2}. Однако тогда должны выполняться соотношения
\begin{flalign*}
	&& g_m(\zeta_m+\beta_m)&<1,&\text{[\eqref{eq:3:3:2}]}&\\
	&& g_m(\zeta_m)&<1,&\text{[\eqref{eq:3:3:2}]}&\\
	&& |(\zeta_m+\beta_m)-\zeta_m|&=\beta_m\\
	&& &<\delta_m,&\text{[\eqref{eq:3:3:1}]}&\\
	&& |f(\zeta_m+\beta_m)-f(\zeta_m)|&=|f_m(\zeta_m+\beta_m)-
		f_m(\zeta_m)|&\text{[\eqref{eq:3:3:3}, \eqref{eq:3:3:2}]}&\\
	&& &=|0-2^{-\mu(m)}|&\text{[\eqref{eq:3:3:3}]}&\\
	&& &=2^{-\mu(m)},
\end{flalign*}
противоречащие сделанным предположениям о свойствах алгорифма \(\mathfrak A_{%
\mu(m)}\).
\end{proof}

\subsection
Заметим, что применительно к многомерному случаю утверждение об осуществимости
функций, интегрируемых по Риману, но не интегрируемых по Дарбу, может быть получено
в качестве следствия из результатов \cite{Dem:1968} о неверности теоремы Фубини
для конструктивного интеграла Римана. Однако такое доказательство требует
привлечения ряда дополнительных представлений, поэтому на его деталях мы здесь
не останавливаемся.

\enlargethispage{\baselineskip}


\begin{thebibliography}{99}
\bibitem{Markov:1974} \emph{А.~А.~Марков.} О языке \(\Ya{\omega |}\)//
ДАН~СССР. "--- 1974. "--- Т.~215, \No~1. "--- С.~57--60.
\bibitem{VD:2009} \emph{А.~А.~Владимиров, М.~Н.~Домбровский--Кабанченко.}
Ступенчатая семантическая система. "--- М.: Изд-во~ВЦ~РАН, 2009.
\bibitem{Ku:1973} \emph{Б.~А.~Кушнер.} Лекции по конструктивному математическому
анализу. "--- М.: Наука, 1973.
\bibitem{ZTs:1962} \emph{И.~Д.~Заславский, Г.~С.~Цейтин.} О сингулярных покрытиях
и связанных с ними свойствах конструктивных функций// Труды~Матем. ин-та
им.~В.~А.~Стеклова. "--- 1962. "--- Т.~67. "--- С.~458--502.
\bibitem{Dem:1973} \emph{О.~Демут.} О конструктивном аналоге связи измеримости
множеств и функций по Лебегу// Comment. Math. Univ.~Carolinae. "--- 1973. "---
B.~14, \No~3. "--- S.~377--396.
\bibitem{Sh:1962} \emph{Н.~А.~Шанин.} Конструктивные вещественные числа
и конструктивные функциональные пространства// Труды~Матем. ин-та
им.~В.~А.~Стеклова. "--- 1962. "--- Т.~67. "--- С.~15--294.
\bibitem{MN:1996} \emph{А.~А.~Марков, Н.~М.~Нагорный.} Теория алгорифмов.
Изд.~2. "--- М.: ФАЗИС, 1996.
\bibitem{Dem:1968} \emph{О.~Демут.} О теореме Фубини для интеграла Римана
в конструктивной математике// Comment. Math. Univ.~Carolinae. "--- 1968. "---
B.~9, \No~4. "--- S.~677--686.
\end{thebibliography}
\end{document}